\documentclass[12pt]{article}
\usepackage{amssymb,latexsym,amsmath}
\usepackage{graphicx}
\def\C{\mathbb C}
\def\R{\mathbb R}
\def\N{\mathbb N}
\def\Z{\mathbb Z}

\newtheorem{thm}{Theorem}[section]
\newtheorem{lem}{Lemma}[section]

\newtheorem{prop}{Proposition}[section]

\begin{document}
\sffamily
\title{Bank-Laine functions with real zeros}
\author{J.K. Langley}
\maketitle

\centerline{\textit{Dedicated to the  memory of Stephan Ruscheweyh}}

\begin{abstract}
Suppose that  $E$ is a real entire function of finite order with  zeros which are all real
but neither bounded above nor bounded below,
such that $E'(z) = \pm 1$ whenever $E(z) = 0$. 
Then either $E$ has an explicit representation in terms of trigonometric functions
or the  zeros of $E$ have exponent of convergence
at least $3$. An example constructed via quasiconformal surgery demonstrates the sharpness of this result.
\end{abstract}

\section{Introduction}

For a non-constant entire function $f$, denote by
$$
\rho (f) = \limsup_{r \to + \infty} \frac{ \log^+ T(r, f)}{\log r} \, , \quad 
\lambda (f) = \limsup_{r \to + \infty} \frac{ \log^+ N(r, 1/f)}{\log r}  \leq \rho (f),
$$
its order of growth and the exponent of convergence of its zeros \cite{Hay2}. 
In their seminal paper \cite{BIL1}, 
Bank and Laine proved several landmark results on the oscillation of solutions of 
\begin{equation}
 \label{de1}
y'' + A(z) y = 0,
\end{equation}
in which $A$ is an entire function. 
Their approach was based on taking linearly independent solutions $f_1, f_2$ of (\ref{de1}), normalised so as to have 
Wronskian $W(f_1, f_2) = f_1 f_2' - f_1' f_2 = 1$, and then considering the product  $E=f_1f_2$, which satisfies
\begin{equation}
 \label{bleq}
 4A = \left( \frac{E'}{E} \right)^2 - 2 \, \frac{E''}{E} - \frac{1}{E^2} . 
\end{equation}
In particular, it was shown in \cite{BIL1} that if $\lambda (E) + \rho (A) < + \infty$ then $\rho(E) < + \infty$, 
whereas if $A$ is transcendental then 
the quotient $U = f_1/f_2$ always has infinite order, since \cite[Chapter 6]{Lai1}
\begin{equation}
 \label{Sfdef}
S_U(z) = \frac{U'''(z)}{U'(z)} - \frac32 \left( \frac{U''(z)}{U'(z)} \right)^2  = 2A,
\end{equation}
where $S_U(z)$ is the Schwarzian derivative. The  following results were proved by Bank and Laine,
Rossi and Shen  \cite{BIL1,Ros,Shen}. 

\begin{thm}[\cite{BIL1,Ros,Shen}]
 \label{thmA}
Let $A$ be an entire function, let $f_1, f_2$ be linearly independent solutions 
of (\ref{de1}) and let $E=f_1f_2$, so that $\lambda (E) = \max \{ \lambda (f_1), \lambda (f_2) \}$.\\
(i) If $A$ is a polynomial of degree $n > 0$ then $\lambda(E) = (n+2)/2 $.\\
(ii) If $\lambda(E) < \rho(A) < + \infty$ 
then $\rho(A) \in \N = \{ 1, 2, \ldots \}$.\\
(iii) If $A$ is transcendental and $\rho(A) \leq 1/2$ then $\lambda(E) = + \infty$, while if $1/2 < \rho(A) < 1$ then 
\begin{equation}
 \label{roineq}
\frac1{\rho(A)} + \frac1{\lambda (E) } \leq 2 .
\end{equation}
\end{thm}
Theorem \ref{thmA}(ii) inspired the 
Bank-Laine conjecture, 
to the effect that if $A$ is a transcendental entire function and $f_1, f_2$ are linearly
independent solutions of $(\ref{de1})$ with $\lambda (f_1f_2)$ finite then  
$\rho (A) \in \N \cup \{ + \infty \} $.
This conjecture has recently been disproved, however, in the first of two remarkable  papers of Bergweiler
and Eremenko 
\cite{bebl1,bebl2} which use  quasiconformal constructions; in the second of these  they show  that  equality is possible in (\ref{roineq}), for every 
choice of $\rho(A) \in (1/2, 1)$.

The  main thrust of this paper concerns the location of zeros of Bank-Laine functions, these being entire functions $E$ such
that $E(z)=0$ implies $E'(z) = \pm 1$. By \cite[Lemma C]{BIL3},
an entire function $E$ is a Bank-Laine function if and only if 
$E=f_1f_2$, where $f_1, f_2$ are linearly independent and solve (\ref{de1}), with $W(f_1, f_2) = 1$  and $A$ entire, satisfying (\ref{bleq}). 
Although 
a Bank-Laine function with unrestricted growth
may have  arbitrary zeros, subject only to these  having no finite limit point \cite{Shen2},
the following is a combination of results from \cite{qcsurg,blnewqc}.

\begin{thm}[\cite{qcsurg,blnewqc}]
 \label{thmDL}
Let $E$ be a Bank-Laine function of finite order, with infinitely many zeros, all real, and denote by $n(r)$ the number of zeros of $E$
in $[-r, r]$.  Then $n(r) \neq o(r)$ as $r \to + \infty$.  If, in addition,  all zeros of $E$ are positive, then $\lambda(E) \geq 3/2$.
\end{thm}

The first assertion of Theorem \ref{thmDL} is evidently sharp, because of $\sin z$, and so is the second, a suitable example 
having been constructed in \cite{blnewqc} using quasiconformal methods. The next theorem will considerably strengthen Theorem \ref{thmDL}
in the case where $E$ is a real Bank-Laine function of finite order with real zeros, these neither bounded above nor bounded below.

\begin{thm}
 \label{thmBLpositive}
Let $E$ be a real Bank-Laine function of finite order, with only real zeros, 
these neither bounded above nor bounded below, 
and let $A$ be the associated coefficient function  in (\ref{de1}) and (\ref{bleq}).
Then one of the following holds:\\
(i) there exist $\eta , \omega_1 , \omega_2 \in \R$ 
such that  $\eta \sin( \omega_1 - \omega_2 ) \neq 0$  and 
\begin{equation}
\label{explicitform}
A = \eta^2, \quad E(z) = \pm  \frac{ \sin (\eta z - \omega_1 ) \sin ( \eta z - \omega_2 ) }
{ \eta \sin ( \omega_1 - \omega_2 )} ;
\end{equation}
(ii) $A$ is transcendental
 and $\lambda(E)  \geq 3$, with $\rho(E) = \rho (A) = 3$  if $\lambda(E) = 3$. 
\end{thm}

Note that Hellerstein, Shen and Williamson \cite{HSW} proved that
 if $A$ is a non-constant polynomial 
then (\ref{de1}) cannot have linearly independent solutions with only real zeros
(see also \cite{Gunreal,Steiradial}).  
A simple example illustrating (i) is given by
$$A = 1/4, \quad E(z) = \sin z = 2 \sin (z/2) \cos (z/2) = -2 \sin(z/2) \sin(z/2 - \pi/2), $$
while 
$E(z) =  \pi^{-1} e^{2 \pi i z^2 } \sin \pi z $  shows that the hypothesis that $E$ is real entire
is not redundant in Theorem \ref{thmBLpositive}. The sharpness of the result will be demonstrated  in Section
\ref{example3}, in which the  
quasiconformal techniques  of \cite{blnewqc} will be adapted
to  construct a real Bank-Laine function $E$,
whose  zeros  are real but  neither bounded above nor  bounded below, such that 
$E$ and the  function $A$   in (\ref{de1}) and (\ref{bleq}) satisfy
$\lambda(E) = \rho(E) = \rho (A) = 3$. 

The author thanks the referees for several helpful suggestions. 

\section{Preliminaries}

Let $G$ be a transcendental meromorphic 
function in the plane and suppose that  $G(z) \to a \in \C \cup \{ \infty \}$ as $z \to \infty$ along a path
$\gamma $; then the inverse $G^{-1}$ is
said to have a transcendental singularity
over the asymptotic value $a$~\cite{BE,Nev}. If $a \in \C$ then for each $\varepsilon > 0$ 
there  exists a component $\Omega = \Omega( a, \varepsilon, G)$ of the
set $\{ z \in \C : |G(z) - a | < \varepsilon \}$ such that
$\gamma \setminus \Omega$ is bounded, and each such $\Omega$ is called a neighbourhood of the singularity \cite{BE}.
Two  paths $\gamma, \gamma'$ on which $G(z) \to a$ determine distinct singularities if the corresponding components
$\Omega( a, \varepsilon, G)$, $\Omega'( a, \varepsilon, G)$ are disjoint for  some $ \varepsilon > 0$.
The singularity 
is called direct  \cite{BE} if $\Omega( a, \varepsilon, G)$, for some
$\varepsilon > 0$, contains no zeros of $G-a$, and logarithmic if 
there exists $\varepsilon > 0$ such that $\log 1/(G-a)$ maps $\Omega( a, \varepsilon, G)$ conformally onto the half-plane 
${\rm Re} \, w > \log 1/\varepsilon$. 
Transcendental singularities over $\infty$ may be 
classified using $1/G$.

Denote by  $B(a, r)$ the open disc of centre $a \in \C$ and radius $r$, and by $Cl(D)$ the closure, with respect to the finite plane,  of   $D \subseteq \C$. 

\begin{prop}
 \label{lemUspeiser}
Let $f_1, f_2$ be  linearly independent solutions of (\ref{de1}), in which $A$ is a transcendental entire function of finite order,
and assume that $W(f_1, f_2) = 1$ and  the zeros of $E = f_1 f_2$ have finite exponent of convergence. Write
\begin{equation}
 \label{Udef1a}
U = \frac{f_2}{f_1}, \quad \frac{U'}{U} = \frac{W(f_1, f_2)}{f_1 f_2} = \frac1E, \quad 
F(z) = \frac{E(z)}{z}.
\end{equation}
Then the following statements hold.\\
(A) $U^{-1}$ has finitely many transcendental singularities over finite non-zero values.\\
(B) $U$ has no critical values and finitely many asymptotic values.\\
(C) Let  $\Omega$ be a  neighbourhood of a
transcendental singularity of $U^{-1}$ over $\alpha \in \C \setminus \{ 0 \}$.
Then $\Omega$
contains a neighbourhood of a direct transcendental singularity of $F^{-1}$ over $\infty$, as well as a path tending to infinity on which $U(z) \to \alpha$
and $F(z) \to \infty$.\\
(D) $U^{-1}$ has infinitely many logarithmic singularities over $0$ or $\infty$.\\
(E) Let $\gamma$ be a path tending to infinity on which $U(z)$ tends to $0$ or $\infty$.
Then  $F(z)$ tends to $0$ on $\gamma$.
\end{prop}
\textit{Proof.} The fact that $U$ has no critical values is well known, and holds since $U'/U = 1/E \neq 0$ and all zeros and poles of $U$ are simple.
Thus (B) follows from  (A), and (A) from \cite[Theorem 1.3]{Lasing2016}, because the Bank-Laine equation (\ref{bleq}) implies that 
$E$ has finite order \cite{BIL1}. 

To prove (C) requires an argument from \cite[Lemma 5.2]{blnewqc}.
Since the singular value $\alpha$ of $U^{-1}$ is isolated, 
the singularity must be logarithmic \cite[p.287]{Nev}.
Hence
there exist $M > 0$ and a component $\Omega_M \subseteq \Omega $ of  
$\{ z \in \C : \, | U(z) - \alpha | < 1/M \}$ which is mapped univalently by $v = \log 1/(U(z)-\alpha)$  onto the half-plane $H_0$ given by ${\rm Re} \, v > \log M$. 
It may be assumed that $M$ is so large that $\Omega_M \cap B(0, 1) = \emptyset$. 
Let $\phi: H_0 \to \Omega_M$ be the inverse function and for $v \in H_0  $ write 
\begin{equation}
 \label{vdef1}
U(z) = \alpha + e^{-v} , 
\quad  E(z) = \frac{U(z)}{U'(z)} = \frac{\alpha  + e^{-v}}{-e^{-v}} \cdot \phi'(v)  = - (1+\alpha e^v) \phi'(v).
\end{equation}
Bieberbach's theorem and 
Koebe's quarter theorem \cite[Chapter 1]{Hay9} give 
$$
\left| \frac{\phi''(v)}{\phi'(v)} \right| \leq \frac4{{\rm Re} \, v - \log M} \quad \hbox{and}
 \quad  \left| \frac{\phi'(v)}{\phi(v)} \right| \leq \frac{4 \pi}{{\rm Re} \, v - \log M} 
$$
on $H_0$, and so there exists 
$c_1 > 0$ such that, as $v \to + \infty$ on $\R$, by (\ref{vdef1}),  
$$
|z| = |\phi(v)| = o( v^{c_1} ) = o( e^v | \phi'(v) | ) = o( |E(z)| ) , \quad \frac{z}{E(z)} = \frac1{F(z)} = o(1) .  
$$
Thus $U(z) \to \alpha$
and $F(z) \to \infty$ as $z \to \infty$ on the image under $\phi$ of the interval $[2 + \log M, + \infty)$. 
On the other hand, if  ${\rm Re} \, v = 1 + \log M $ then, again by (\ref{vdef1}),
$$
|F(z)| \leq  ( 1 + | \alpha | M e )  \left| \frac{\phi'(v)}{\phi(v)} \right| \leq  ( 1 + | \alpha | M e ) 4 \pi .
$$
Hence there exist large positive $M_0, M_1 $ and a component $C_0$ of 
$\{ v \in H_0: \, |F(\phi(v))| > M_0 \}$ containing an interval $[M_1, + \infty)$, such that $Cl(C_0) \subseteq H_0$, and
$\phi(C_0)$ is the required neighbourhood  of a direct  singularity of $F^{-1}$ over $\infty$.

Next, (D) follows from (A) and the result of Nevanlinna-Elfving \cite{elf,Nev2}, which implies that if $U^{-1}$ 
has finitely many transcendental singularities then its Schwarzian derivative $2A$ must be a rational function, contrary to hypothesis
(see also \cite{schwarzian}). 

Finally, to prove (E), take a path $\gamma$ tending to infinity on which $V(z) \to \infty$, where $V$ is $U$ or $1/U$. 
Since $\infty$ is not a limit point of singular values of $V^{-1}$, a standard estimate (see \cite{EL}
or \cite[Section 6]{sixsmithEL}) gives positive constants $c_2, c_3$
such that,  as $z $ tends to infinity on $\gamma$, 
$$
\frac1{|F(z)|} = 
\left| \frac{z}{E(z)} \right| = \left| \frac{zV'(z)}{V(z)} \right| \geq
c_2 \log \left| \frac{V(z)}{c_3} \right| \to + \infty . 
$$
\hfill$\Box$
\vspace{.1in}

\section{Proof of Theorem \ref{thmBLpositive}}

Suppose that $E$ and $A$ are as in the hypotheses. Then there  exist  solutions $f_1, f_2$ of (\ref{de1}) such that 
$W(f_1, f_2) = 1$ and $ E = f_1 f_2$.  Furthermore, $A \not \equiv 0$ and
each $f_j$ has infinitely many zeros on each of the positive and negative real axes, 
since $f_j(z) = 0$ gives $E'(z) = (-1)^j$ and 
the sign of $E'$ at successive zeros must alternate.

If  $A = \eta^2$ is constant,
then each $f_j$ is a linear combination of $e^{i \eta z} $ and $e^{- i \eta z }$;
moreover, since the $f_j$ have infinitely many real zeros,
 $\eta $ must be real,  and 
$f_j(z) = A_j \sin ( \eta z - \omega_j ) $, with $A_j , \omega_j$ constants and $\omega_j$ real,
which gives $E(\omega_1/\eta ) = 0$ and forces  $E'(\omega_1/ \eta ) = \pm 1$ and (\ref{explicitform}).
Assume henceforth that $A$ is non-constant.


\begin{lem}
 \label{lemAtrans}
The  function $A$     is transcendental.
\end{lem}
\textit{Proof.}
This follows from results in \cite{Gunreal,HSW, Steiradial}, and may be proved via
the following slight modification of  \cite[Lemma 5.1]{blnewqc}.
Suppose that $A$ is a polynomial in (\ref{de1}), non-constant by assumption,
satisfying $A(z) = a_n z^n (1+o(1))$ as $z \to \infty$. Then  
there are $n+2 > 2$  critical rays given by $\arg z = \theta^*$, where $a_n e^{i (n+2) \theta^*} $ is real and positive, and a combination of
the Liouville transformation
\begin{equation*}
 \label{ltdef}
Y(Z) = A(z)^{1/4} y(z) , 
\quad Z = \int^z \, A(t)^{1/2} \, dt ,
\end{equation*}
with Hille's asymptotic method \cite{Hil2}  generates linearly independent \textit{principal} solutions of (\ref{de1}) given by
$A(z)^{-1/4} e^{\pm iZ} (1+o(1))$ on  sectors symmetric about these rays.
On one side of the critical ray, one  of these  principal solutions
is large, while the other is small, 
these roles being reversed as the ray is crossed. 
Since the $f_j$ have infinitely many positive zeros, the positive real axis must be one of these $n + 2$  critical rays,
and each $f_j$ must be a non-trivial linear combination of the two principal solutions and so
large in both  adjacent sectors.  
Let $L$ be the first other critical ray encountered on
moving counter-clockwise from the positive real axis. Then $L$ is not the negative real axis, as $n > 0$.
Since the $f_j$ have only real zeros, both  must change from large to  small as  $L$ is crossed.
But $f_1, f_2$ cannot be small in the same sector, because $W(f_1, f_2) = 1$.
\hfill$\Box$
\vspace{.1in}

Assume henceforth that $A$ is transcendental, but that $\lambda (E) =
\lambda \leq 3$. Then   the canonical product $\Pi_0$ over the zeros of $E$ has order  $\lambda$, and  there exists a  real polynomial $P_0$ with
$E = \Pi_0 \exp( P_0 )$.
If $P_0$ has degree greater than $\lambda$,  then 
the zeros of $E$ have Nevanlinna deficiency $\delta(0, E) = 1$ \cite[p.42]{Hay2}, contradicting \cite[Theorem 4.1]{qcsurg} 
(see also \cite[Theorem 2.1]{Lasparse}). Assume henceforth that $E$ 
has order  $\rho (E) = \lambda \leq 3$: then 
$\rho(A) \leq \lambda  \leq 3$ by (\ref{bleq}).  

Define $U$ and $F$ by~(\ref{Udef1a}). 
Since $U'/U = 1/E$ is real, there exists $\theta \in \R$ such that $U = f_2/f_1 = e^{ 2 i \theta } U_0$, with $U_0$ real meromorphic.
Replacing $f_1 $ by $f_1 e^{i \theta}$ and $f_2 $ by $f_2 e^{- i \theta}$ leaves $E$ unchanged; hence it may be assumed that $\theta=0$ and $U$ is real meromorphic.
By Proposition \ref{lemUspeiser}, $U^{-1}$ has finitely many transcendental singularities over finite non-zero values,
but infinitely many transcendental singularities over $0$ or $\infty$. 

\begin{lem}
 \label{lemsingU}
$U^{-1}$ has at least four  logarithmic singularities over finite non-real values. 
\end{lem}
\textit{Proof.}
Note that it is not asserted that the corresponding four asymptotic values must all be distinct.
Take  zeros $x_0 , x_1, x_2 \in \R$ of $f_2$, with $0 < x_0 < x_1 < x_2$,
and  the supremum $R$ of all $r > 0$ such that the branch of $U^{-1}$ mapping $0$ to $x_1$ extends analytically to $B(0, r)$. Then $R < + \infty$ and 
$U$ maps a simply connected domain $\Omega_1$, with $x_1 \in \Omega_1$, univalently onto
$B(0, R)$. Moreover, $U^{-1}$ has a singularity over some $\alpha$ with $|\alpha | = R$, and so 
$\Omega_1$ contains a path $\gamma$ which  tends to infinity,
mapped by $U$ onto the half-open line segment $[0, \alpha)$. If $\alpha \in \R$ then, since $U$ is real meromorphic and univalent on $\Omega_1$, 
the path  $\gamma $ must be  $(-\infty, x_1]$ or $[x_1, + \infty)$, a 
contradiction since  $x_0, x_2 \not \in \gamma$. 
Hence 
$\alpha \not \in \R$ and so  $U^{-1}$ has logarithmic singularities over $\alpha$ and $\overline{\alpha}$, by Proposition \ref{lemUspeiser}(B)
and  \cite[p.287]{Nev}.

Suppose now that $U^{-1}$ has no other  logarithmic singularities over finite non-real values. Then, without loss of generality, there exist 
neighbourhoods $\Omega_\alpha \subseteq H^+$ and $\Omega_{\overline{\alpha}}  \subseteq H^-$ of the  singularities over $\alpha$ and $\overline{\alpha}$
respectively, where $H^+, H^-$ denote the upper and lower half-planes.
The argument of the previous paragraph shows that all but finitely many zeros of $f_2$ are joined to $\Omega_\alpha $
by a path which is mapped by $U$ onto  $[0, \alpha)$, and to $\Omega_{\overline{\alpha}} $ by a path 
 mapped  onto  $[0, \overline{\alpha})$. Since the set of zeros of $f_2$ is neither bounded above nor bounded below,
this excludes transcendental singularities of $U^{-1}$ over $\infty$, and an almost identical argument applied to $f_1$ 
rules out transcendental singularities of $U^{-1}$ over $0$. This contradiction proves the lemma. 
\hfill$\Box$
\vspace{.1in}

Thus, without loss of generality, there exist disjoint 
neighbourhoods $U_1, U_2 \subseteq H^+$ and $U_3, U_4 \subseteq H^-$ 
of  singularities of $U^{-1}$ over  values $\alpha_j \in \C \setminus \R$, $j = 1, \ldots, 4$.
Proposition \ref{lemUspeiser} delivers for each $j$  a path $\tau_j \subseteq U_j$ on which
$U(z) \to \alpha_j$ and $F(z) \to \infty$. Take a circle $|z| = R$, with $R$ large, 
which meets all four  $\tau_j$ and on which $F$ has no zeros. Thus $|F(z)|$ is bounded below on the union of the circle and the $\tau_j$.
Since $\tau_1, \tau_2$  lie in neighbourhoods of distinct singularities, while $U$ has only real poles,
there must exist a path tending to infinity in $|z| > R$, lying between 
$\tau_1 $ and $\tau_2$, on which $U(z) \to \infty$.
Because $F$ has only real zeros, Proposition \ref{lemUspeiser} now gives neighbourhoods $V_j \subseteq U_j$, for $j = 1, \ldots, 4$, 
of  direct singularities of $F^{-1}$ over $\infty$, 
and neighbourhoods $V_5 \subseteq H^+$ and $V_6 \subseteq H^-$ of
direct singularities of
$F^{-1}$ over  $0$.

This gives positive constants $M_j$ and non-constant, non-negative subharmonic functions $u_1, \ldots, u_6$, with pairwise disjoint supports $V_j$, 
such that $u_j = \log |F/M_j|$ on $V_j$, for $j=1, \ldots, 4$,
while $u_j = \log |M_j/F|$ on $V_j$, for $j = 5, 6$. 
Thus 
$u_1, \ldots, u_4$ have order  $\rho (u_j) \leq \rho(F) \leq \rho(E) = \lambda(E) \leq
3$. Moreover, $u_5, u_6$ have order $\rho (u_j) \leq \rho(A) \leq  3$, because (\ref{bleq}) and Poisson's formula yield as $r \to + \infty$, for $k = 5, 6$,
\begin{eqnarray*}
 \max \{ u_k(z): \, |z| = r \} &\leq& \frac3{2 \pi} \int_0^{2 \pi} u_k(2r e^{i t} ) \, dt \leq  3 m(2r, 1/F) + O(1)  \\
&\leq& 3 m(2r, 1/E) + O( \log r) \leq 3 T(2r, A) + O( \log r) .
\end{eqnarray*}
For $j = 1, \ldots, 6$ and $t > 0$ let $\theta_j(t)$ be the angular measure 
of  $\{ z \in \C : \, |z| = t, \, u_j(t) > 0 \}$.
Let $S$ be large and positive: then
a well known consequence of Carleman's 
estimate for harmonic measure \cite[pp.116-7]{Tsuji} gives, as $r \to + \infty$,
\begin{eqnarray*}
 36 \log \frac{r}{S} &= & \int_S^r   \left( \sum_{j=1}^6 1 \right)^2 \, \frac{dt}t   \leq
 \int_S^r  \left( \sum_{j=1}^6 \theta_j(t) \right) \left( \sum_{j=1}^6 \frac1{\theta_j(t)} \right) \, \frac{dt}t  \\
&\leq& 2 \sum_{j=1}^6 \int_S^r  \frac\pi{t \theta_j(t)}  \, dt    \leq 2 \sum_{j=1}^6 \log ( \max \{ u_j(z): \, |z| = 2r \}  )  + O(1) 
\\
&\leq& 2   \sum_{j=1}^6 (\rho (u_j) + o(1)) \log r \leq ( 8 \lambda (E) + 4 \rho (A) + o(1)) \log r \leq ( 36 + o(1)) \log r .
\end{eqnarray*}
It follows at once that 
$\lambda (E) = \rho (E) = \rho (A) = 3$. 
\hfill$\Box$
\vspace{.1in}


\section{A real Bank-Laine function  with real zeros}\label{example3}

The construction of an example demonstrating that Theorem \ref{thmBLpositive} is sharp
starts with the following. 

\begin{lem}
 \label{lem20a}
The  M\"obius transformation
\begin{equation}
 \label{301}
w = T(v) = \frac{e^{i\pi/4}(1+e^{i\pi/4}v)}{1-e^{-i\pi/4}v} = \frac{e^{-i\pi/4}(v - e^{i3 \pi/4})}{v-e^{i\pi/4}} 
\end{equation}
satisfies
\begin{equation}
\label{301a}
 T( e^{i \pi/4} )  = \infty, \quad T(e^{i3\pi/4}) = 0, \quad T(i) = -1, \quad T(0) = e^{i \pi/4}, \quad T(\infty) = e^{-i \pi/4}.
\end{equation}
In addition, 
$T$ maps the unit circle $|v| = 1$ onto $\R \cup \{ \infty \}$, 
with  ${\rm Im} \, T(v) > 0$ for $|v| < 1$. 
Moreover, ${\rm Re} \, v = 0$ implies that $|T(v)| = 1$, while $|T(v)| > 1$ for ${\rm Re} \, v > 0$, and 
$T$ maps the line segment $[0, i]$ onto the counter-clockwise circular arc from $e^{i \pi/4}$ to $-1$.
\end{lem}
\textit{Proof.} All assertions  follow from (\ref{301a}) and the following observations: first,  
(\ref{301}) implies that $|T(v)| > 1$ precisely when 
$v $ is further from $e^{i 3 \pi/4}$ than from $e^{i  \pi/4}$; second, as $v$ describes the positive imaginary axis, 
$w$ travels around the unit circle from $e^{i \pi/4}$ to $e^{-i \pi/4}$ via $-1$.
\hfill$\Box$
\vspace{.1in}

\begin{lem}
 \label{lem20b}
Write $u = s+it$ with $ s, t \in \R$.
Then the   locally univalent functions 
\begin{equation}
  \label{201}
f_1(u) = e^{i \pi/4} \exp(  \sqrt{2} \, v ),  \quad 
f_2(u) = T(v) , \quad v = e^{iu} ,
\end{equation}
have the following properties:\\
(A)  ${\rm Im} \, f_2(u) > 0$ for $t > 0$, and $f_2$  is a piecewise increasing mapping from $\R$ to $\R \cup \{ \infty \}$;\\
(B) $f_2$ has asymptotic values $e^{\pm i \pi/4}$ and poles at $u = (2k+1/4) \pi $, $k \in \Z$,
as well as zeros at $u = (2k+3/4) \pi$, $k \in \Z$;\\
(C) for  $j = 1, 2$,  the function
$\log |f_j(u)| $ is  positive for $- \pi/2 < s <  \pi /2$, and negative if $-3 \pi/2 < s < - \pi/2$ or $ \pi/2 < s < 3 \pi/2$;\\
(D) $f_1$ maps the vertical line $\gamma_1$ given by $s = -\pi $ onto the open line segment $(0, e^{i \pi/4} )$ and
there exists a path $\gamma_2$, starting at $3 \pi/4$ and tending to infinity in the half-strip $\pi /2 < s < 3 \pi/2$, $0 \leq t < + \infty$, 
which is mapped by $f_2$ onto the half-open line segment $[0, e^{i \pi/4} )$.
\end{lem}
\textit{Proof.} (A) and (B) follow from Lemma \ref{lem20a}, as does (C) for $f_2$, while (C) for $f_1$ is an immediate consequence of the formula
$\log |f_1(u)| = \sqrt{2} e^{-t} \cos s$. 
The assertion (D) for $f_1$ is obvious, while (D) for $f_2$ follows from (A) and 
analytic continuation of $f_2^{-1}$, the  only
singular values  of which are the two values omitted by $f_2$, namely $e^{\pm i \pi/4}$.
\hfill$\Box$
\vspace{.1in} 

The construction will 
proceed by first forming, on the sector $0 < \arg u < 3 \pi /2$, a quasiregular mapping
which is $f_1(u) $ for ${\rm Re} \, u \leq - \pi/2$ and
$f_2(u) $ for ${\rm Re} \, u \geq \pi/2$. A modification of this mapping will
be pulled back to the first quadrant, which will then
permit extension via double reflection to a quasimeromorphic mapping
on the whole plane.  Application of 
the Teichm\"uller-Belinskii theorem \cite{LV} will result in a locally univalent meromorphic
function $U$ for which $E= U/U'$ will be the required Bank-Laine function.
To this end, set 
\begin{eqnarray}
 \label{D1def}
D_0 &=& \{ u \in \C : \, 0 < |u| < + \infty , \, 0 < \arg u < 3 \pi /2 \} , \nonumber \\
D_1 &=&  
E_1 \cup E_2 ,\nonumber \\
E_1 &=&  \{  s + it: \, - \pi /2 < s < 0, \, - \infty < t < + \infty  \} , \nonumber \\
E_2 &=& \{  s + it: \, - \pi /2 < s < \pi/2 , \, 0 < t < + \infty  \} , \nonumber \\
D_2 &=& \{ v \in \C : \, 0 < |v| < + \infty, \, - \pi/2 < \arg v < 0 \} ,\nonumber \\
D_3 &=& D_2 \cup \{ \zeta \in \C : \, | \zeta | < 1, \, {\rm Re} \, \zeta > 0 \}, \nonumber \\
D_4 &=& \{ \sigma + i \tau : \, 0 < \sigma < + \infty, \,  - \infty < \tau <  \pi \}.
\end{eqnarray}
The following is \cite[Lemma 6.1]{blnewqc}.


\begin{lem}[\cite{blnewqc}]
 \label{lemqc2} 
Let $h: (-\infty, 1] \to (-\infty, 0] $ be a continuous bijection, such that $h(1)=0$ while $h'$ is continuous and has positive upper and lower bounds
for $- \infty < y < 1$ (that is, there exists 
$\varepsilon > 0$ such that $\varepsilon < h'(y) < 1/\varepsilon$ for $- \infty < y < 1$). Then there exists a homeomorphism $\psi $ from the closure of $D_3$ to that of 
$D_2$, such that: (A) $\psi$ maps $D_3$ quasiconformally onto $D_2$, with $\psi(z) \to \infty$ and $\psi(z) = O( |z| )$ as $z \to \infty$ in $D_3$; (B)
$\psi(iy) = i h(y) $ for $-\infty < y \leq 1$; (C) $\psi(z)$ is real and strictly increasing as $z  $ describes  
the boundary of $D_3$ clockwise from  $i$ to infinity. 
\end{lem}
\hfill$\Box$
\vspace{.1in}

\begin{lem}
 \label{lemqrmap}
With $D_0, D_1$  as in (\ref{D1def}), and the  $f_j$ as in Lemma \ref{lem20b}, let $E_0= Cl(D_0)$ and define $F$ on $ E_0 \setminus D_1 $ by 
\begin{eqnarray}
 \label{204}
F(s +it) &=& f_1(s +it) \quad \text{for $- \infty < s \leq -\pi/2$,  $t \in \R$}, \nonumber \\
F(s +it) &=& f_2(s +it) \quad \text{for $\pi/2 \leq s < + \infty $,  $0 \leq t < + \infty $}.
\end{eqnarray}
Then $F$ extends to a  mapping from $E_0$ into the extended plane, continuous with respect to the spherical metric,
with the following properties.\\
(i) $H = \log F$ maps $D_1$ quasiconformally onto  $D_4$, with $H(\pi/2) = i \pi $.\\
(ii) $F$ is locally injective on $E_0$. \\
(iii) 
Let $L_0$ be the path consisting of the line segment from $3\pi /4 $ to $0$ followed by the negative imaginary axis in the direction of
$- i \infty$. Then $F(3 \pi/4) = 0$ and $F(u)$ is real and strictly decreasing as $u$ describes $L_0$, mapping 
$L_0$ onto the non-positive real axis. Moreover, each $u_0 \in L_0$ has $s_0 > 0$ such that 
${\rm Im} \, F(u) > 0$ on $D_0 \cap B(u_0, s_0)$.
\\
(iv) There exists $c > 0$ such that $|F(u)| \leq \exp \exp ( c |u| )$ for $u \in D_0$ lying on the circles $|u| = n \pi $, $n  \in \N$. 
\end{lem}
\textit{Proof.} 
First, observe that $v = e^{iu}$ maps $D_1$ onto $D_3$, with $v( \pi/2) = i$ and $v \to 0$ as
${\rm Im} \, u \to + \infty$, as well as $v \to \infty$ as $ {\rm Im} \, u \to - \infty$. Indeed,
the boundary of $D_1$  is mapped by $v = e^{iu}$ as follows: the line ${\rm Re} \, u = - \pi/2$  to the negative imaginary axis;
the half-line ${\rm Re} \, u =  \pi/2$, $0 \leq {\rm Im} \, u < + \infty$,  to the segment 
$v=iy$, $0 < y \leq 1$; the  interval $[0, \pi/2] \subseteq \R$  to the arc of the unit circle from $1$ to $i$; the negative imaginary axis  to
$(1, + \infty)$. Using Lemma \ref{lem20a} and the principal argument, set 
\begin{equation}
\label{hmatch}
h(y) = 
 \begin{cases}
 \frac{\pi}4 + \sqrt{2} \, y  &\text{for $- \infty < y \leq 0$,}\\
 \arg T(iy) =  \frac{\pi}4 + \arg \left(  \frac{1+e^{i\pi/4}i y }{1-e^{-i\pi/4} iy } \right) =
\frac{\pi}4   - i \log \left(  \frac{1+e^{i3 \pi/4} y }{1 + e^{-i3 \pi/4} y } \right)   &\text{for $0 < y \leq 1$.}
\end{cases}
\end{equation}
Then $h(1) = \pi $ and, for $0 < y < 1$, 
$$
 h'(y) = - i  \left(  \frac{e^{i3 \pi/4} }{1+e^{i3 \pi/4} y }-  \frac{e^{- i3 \pi/4} }{1 + e^{-i3 \pi/4} y } \right) = \frac{ 2  \sin (3 \pi/4) }{| 1+e^{i3 \pi/4} y |^2 } > 0,
$$
so that 
$\lim_{y \to 1-} h'(y)$  is finite but positive, and  $\lim_{y \to 0+} h'(y) = \sqrt{2}$, which leads to $h'(0) = \sqrt{2}$. 
Thus $h$ is  a continuous bijection from $(-\infty, 1]$ to $(-\infty, \pi ] $ and
$h'$ exists and is continuous on  $( -\infty , 1)$, with positive upper and lower bounds there. Applying Lemma \ref{lemqc2} to $h(y) - \pi $
gives 
a homeomorphism $\psi $ from the closure of $D_3$ to that of the 
quadrant $D_4$ in (\ref{D1def}), 
such that $\psi$ maps $D_3$ quasiconformally onto $D_4$, 
with $\psi(v) = O( |v| )$ as $v \to \infty$ in $D_3$ and $\psi(iy) = i h(y) $ for $-\infty < y \leq 1$. 
The function $G = \exp \circ \psi $ is then continuous  on $Cl(D_3)$ and satisfies, by (\ref{hmatch}),  
\begin{eqnarray}
G(v)  &=&  \exp ( i h(y) ) = e^{i \pi/4}  \exp( \sqrt{2} \, iy  ) = e^{i \pi/4}  \exp( \sqrt{2}  \, v ) \quad \text{for $v = iy$, $- \infty < y \leq 0$}, \nonumber \\
G(v)  &=& \exp ( i h(y) ) =  T(iy) = T(v)  \quad \text{for $v=iy$, $0 < y \leq  1$}.
\label{Gmatch}
\end{eqnarray}


Now set $F(u) = G(e^{iu}) $ on $Cl(D_1)$. Then (\ref{301}), (\ref{201}), (\ref{204}),
(\ref{Gmatch}) and the  properties already noted of the mapping $v = e^{iu}$ 
from $D_1$ to $D_3$ ensure that
$F$ is well-defined and continuous
on  $E_0$, and  that (i) holds.
Because $\psi$ is injective on $D_3$, and $|F| > 1$ on $D_1$, Lemma \ref{lem20b}(C)
implies  (ii). 
To establish (iii), observe first that $F(u) = f_2(u)$ is real and decreases from $0$ to $-1$ 
as $u$ traverses the line segment from $3 \pi/4$ to $\pi/2$, by (\ref{301a}), (\ref{201}) and
Lemma \ref{lem20b}(A). Next, as $u$ follows $L_0$ from $ \pi/2$ towards infinity, $v = e^{iu}$ 
describes the boundary of $D_3$ clockwise from  $i$ to infinity, so that $\psi(v)$ 
describes the half-line $\{ \sigma + i \pi : \, 0 \leq  \sigma < + \infty \}$, by Lemma \ref{lemqc2} applied to  $h(y) - \pi $, 
and $F(u) = G(v) = \exp ( \psi (v) )$ travels from $-1$ along the negative real axis towards $- \infty$.
Finally, to prove (iv), note first that (\ref{204}) and Lemma \ref{lem20b} show that it is enough to bound  $|F(u)|$ for $u \in D_1$, 
and hence it suffices to consider $G(v)$, which  is continuous  on $Cl(D_3)$ and satisfies,
as $v = e^{iu} \to \infty$ in  $D_3$,
$$
|F(u)| = |G(v)|\leq \exp( |\psi(v)| )  \leq \exp( O(|v| ) ) = \exp \left( O( |e^{iu} | )  \right) \leq \exp \exp ( 2 |u| ).
$$
\hfill$\Box$
\vspace{.1in} 

Next, let $L$ be the  quasiconformal mapping of the extended plane given by $L(re^{i \theta}) = r e^{i g(\theta)}$ for $r > 0$ and $0 \leq \theta \leq 2 \pi$,
where $g$ is  continuous, strictly increasing and piecewise linear with
$$
g(\theta) = \theta \quad \hbox{for} \quad 0 \leq \theta \leq \pi/3,  
\quad g(\pi) = \pi/2, \quad  g(2 \pi ) = 2 \pi  .
$$
Let $E_3$ be the component of 
$E_0 \setminus ( \gamma_1 \cup \gamma_2 )$ which contains $D_1$, where $\gamma_1, \gamma_2$ are as in Lemma \ref{lem20b}(D),  and set 
$V(u) = L(F(u))$ for  $u \in E_3$, with $V(u) = F(u)$ on $ E_0 \setminus E_3$.
Since $F$ maps $\gamma_1 \cup \gamma_2 $ into the segment $[0, e^{i \pi/4})$, on which
$L$ is the identity, $V$ is well-defined and continuous on $E_0$,
and quasiregular and non-zero on $D_0$, with $V(u) = F(u) = f_2(u) \in \R \cup \{ \infty \}$  on $[ 3 \pi/4, + \infty)$. Furthermore,
$V$ maps  the path $L_0$ in Lemma \ref{lemqrmap}(iii) onto the non-negative imaginary axis, each $u_0 \in L_0$ having $s_0 > 0$
such that $0 < \arg V(u) < \pi /2$ on $D_0 \cap B(u_0, s_0)$.

Set  $x_0 = (3 \pi/4)^{2/3} \in (0, + \infty)$, 
and on  the  quadrant $D_5$ given by $0 < \arg z < \pi/2$
write  $\zeta = x_0 + z^2$ and   $u = \eta (z) = \zeta^{3/2}$, taking the principal branch. Then $\eta$ maps $D_5$ onto $D_0$ 
and extends continuously to $\partial D_5$, with $0$ mapped to $ 3 \pi/4$, the non-negative imaginary axis to
the path $L_0$, and the non-negative real axis to $[ 3 \pi/4, + \infty)$. 
Set $Y(z) = V(\eta(z)) = V(( x_0 + z^2)^{3/2} ) $ 
on $D_5$ and extend $Y$ to $Cl(D_5)$  by continuity. Then $Y$ 
maps the non-negative real axis into  $\R \cup \{ \infty \}$, with $Y(0) = V(3\pi/4) = 0$,
and is a bijection from  the non-negative imaginary axis to itself.
Double reflection, first across the imaginary axis and then across the real axis,   extends $Y$  to the whole plane.
The resulting function is  locally injective in the plane, by Lemma \ref{lem20b}(A) and the mapping properties of   $V$ on $L_0$, and
quasimeromorphic  \cite[Ch. I, Theorem 8.3]{LV}. Further, $Y$ now maps $\R$ into $\R \cup \{ \infty \}$, and has only real zeros and poles.
If $x $ is large and positive then $x$ is a zero or pole of $Y$ if and only if 
$\eta (x) \sim x^3 $ is a zero or pole of $f_2$. Thus the set of zeros and poles of $Y$ is neither bounded above nor bounded below,
and by Lemma \ref{lem20b}  the number $n_Y(r)$ of these  
in $[-r, r]$ satisfies 
\begin{equation}
 \label{Vzerosest}
c_1 r^3 \leq n_Y(r) \leq c_2 r^3 \quad \hbox{as} \quad r \to + \infty ,
\end{equation}
in which the $c_j$ denote positive constants. 
Moreover,  Lemma \ref{lemqrmap}(iv) gives, for large $n \in \N$,
\begin{equation}
 \label{rndef}
\log^+ \log^+ |Y(z)| =  O( n ) ,
\end{equation}
initially for $ z \in D_5$ with $| x_0 + z^2 | =  \left(n \pi  \right)^{2/3} $, and hence by reflection on a Jordan curve $\Gamma_n$ 
on which $|z| \sim  \left(n \pi  \right)^{1/3}$.

The remainder of the construction proceeds as in \cite{bebl1,blnewqc}. 
Let $E_4$ be the pre-image in  $D_5$ of $D_6 = \{ u \in D_0 : \, -2 \pi < {\rm Re} \, u < 2 \pi \}$
under $u = \eta (z)$.
If $E_4'= \{ z \in E_4 : |z| > R' \}$, where $R'$ is large, then writing $z = x+iy$,
$u = \kappa + i \lambda $, with $x, y, \kappa , \lambda $ real, leads to
\begin{eqnarray}
 \label{teich1}
\int_{E_4'} \frac1{|z|^2} \, dx dy &=& 
\int_{\eta\left(E_4'\right)} \frac1{|z \eta'(z)|^2 } \, d\kappa d \lambda =  \int_{\eta\left(E_4'\right)} \frac1{ | 9 z^4 (x_0+z^2)|} \, d\kappa d \lambda \nonumber \\
&\leq& c_3 + c_4 \int_{u \in D_6, |u| > 1} \frac1{|u|^2} \, d\kappa d \lambda < + \infty . 
\end{eqnarray}
Now let $F_4$ be the closure of the union of $E_4$ and its reflections across 
the  real and imaginary axes. 
Then $Y$ is meromorphic off $F_4$ and
(\ref{teich1}) implies that
the complex dilatation $\mu_Y$ of $Y$ satisfies
\begin{equation}
 \int_{1 \leq |z| < + \infty} \left| \frac{\mu_Y(z)}{z^2} \right| \, dx dy  \leq \int_{1 \leq |z| < + \infty, z \in F_4} \frac{1}{|z|^2} \, dx dy   < + \infty .
 \label{teich}
\end{equation}
Let $\phi$ be the unique quasiconformal homeomorphism of the extended plane which solves the Beltrami equation 
$\phi_{\overline{z}} =  \mu_Y \phi_z $ a.e.
and fixes each of $0$, $1$ and $\infty$ \cite{LV}. In view of (\ref{teich}) and the Teichm\"uller-Belinskii theorem \cite[Ch. V, Theorem 6.1]{LV},
there exists $\alpha\in \C \setminus \{ 0 \}$ with 
\begin{equation}
 \label{phiasym}
\phi(z) \sim \alpha z
\end{equation}
as $z \to \infty$. Furthermore, there exists a 
locally univalent meromorphic function $U$ such that $Y = U \circ \phi$ on $\C$, and writing  $U_1(z) = \overline{U(\overline{z})}$ gives
$$
U(\phi(z)) = Y(z) = \overline{Y(\overline{z})} =  \overline{U(\phi(\overline{z}))} = U_1 \left( \overline{\phi(\overline{z})} \right) .
$$
Thus $\phi(z)$ and $\overline{\phi(\overline{z})}$ have the same complex dilation a.e. and, since  both fix $0$, $1$ and $\infty$,
they must agree, so that $\phi$ is real on $\R$ and $U$ is real meromorphic. Moreover, 
all zeros and poles of $U$ are real, and $E = U/U'$ is a real Bank-Laine
function with real zeros, these neither bounded above nor bounded below. 
Let $\Pi_1$ and $\Pi_2$ be the canonical products over the zeros and poles of $U$ respectively. 
Then (\ref{Vzerosest}) and (\ref{phiasym}) imply that $\Pi_1$ and $\Pi_2$ have order at most $3$, and 
that the associated coefficient function $A$ in (\ref{de1}) and (\ref{bleq}) cannot be constant. 
There exists  an entire function $Q$ such that 
\begin{equation}
 \label{hdef}
U = \frac{\Pi_1}{\Pi_2} \, e^Q, \quad 
\frac1E = \frac{U'}{U} = \frac{\Pi_1'}{\Pi_1} -  \frac{\Pi_2'}{\Pi_2} + Q' .
\end{equation}
By (\ref{rndef}) and (\ref{phiasym}), the entire function $f_0 = \Pi_2 U$ satisfies, on the Jordan curve $\Gamma_n$,  
$$
| f_0 ( \phi(z) )| = | \Pi_2( \phi (z) ) Y(z) | \leq C_0 e^{|\phi(z)|^4} \exp \exp ( C_1 n )  \leq \exp \exp ( C_2 n ),
$$
in  which the  positive constants $C_j $ are independent of $n$, and so $| f_0 ( w )| \leq \exp \exp ( C_2 n )$
for $w$ on $\phi(\Gamma_n)$. Since $\phi(\Gamma_n)$ encloses a circle $|w| = C_3 n^{1/3} $,  the maximum principle gives
\begin{equation*}
\log T(r,f_0 ) \leq \log \log M(r,f_0 ) = O(r^{3}) \quad \text{as $r \to + \infty$}. 
\end{equation*}
On combination with (\ref{hdef}) and the lemma of the logarithmic derivative \cite{Hay2}, this leads to 
$$T(r, Q') = m(r, Q') \leq m \left(r, \frac{f_0'}{f_0 } \right) + m \left(r, \frac{\Pi_1'}{\Pi_1} \right) + O(1) 
= O( r^{3} ) \quad \text{as $r \to + \infty$}.
$$
Hence (\ref{hdef}) implies that $\rho(E) \leq 3$, and $\lambda (E) = \rho (E) = \rho (A) = 3$ by Theorem \ref{thmBLpositive}. 
\hfill$\Box$

\noindent
School of Mathematical Sciences, University of Nottingham, NG7 2RD.\\
james.langley@nottingham.ac.uk


\vspace{.1in}



{\footnotesize
\begin{thebibliography}{99}
\bibitem{BIL1}
S. Bank and I. Laine, On the oscillation theory of
$f''	  + Af = 0$ where $A$ is entire, Trans. Amer. Math. Soc. 273
(1982), 351-363.
\bibitem{BIL3} S. Bank and I. Laine, On the zeros of meromorphic solutions of second-order linear differential equations, Comment. Math. Helv. 58 (1983), 656-677.
\bibitem{BE}
W. Bergweiler and A. Eremenko, On the singularities of the inverse to
a meromorphic function of finite order, Rev. Mat. Iberoamericana 11 (1995), 355-373.
\bibitem{bebl1}
W. Bergweiler and A. Eremenko, On the Bank-Laine conjecture,  J. Eur. Math. Soc. 19 (2017), 1899-1909.
\bibitem{bebl2}
W. Bergweiler and A. Eremenko, Quasiconformal surgery and linear differential equations,  J. Analyse Math. 137 (2019), 751-812.
\bibitem{qcsurg}
D. Drasin and J.K. Langley, Bank-Laine functions via quasiconformal surgery, Transcendental Dynamics and Complex Analysis, London 
Math. Soc. Lecture Notes 348 (2008), Cambridge University Press, 165-178. 
\bibitem{elf}
G. Elfving, \"Uber eine Klasse von Riemannschen Fl\"achen und ihre Uniformisierung, Acta Soc. Sci. Fenn. 2 (1934) 1-60.
\bibitem{EL} A. E. Eremenko and M.Yu. Lyubich, Dynamical properties of some classes of entire functions, Ann. Inst. Fourier Grenoble 42 (1992), 989-1020.
\bibitem{Gunreal}
G. G. Gundersen, On the real zeros of solutions of $f'' + A(z) f = 0$ where $f$ is entire,
Ann. Acad. Sci. Fenn. Ser. A. I. Math. 11 (1986), 275-294. 
(2) 37 (1988), 88-104.
\bibitem{Hay2}
W.K. Hayman, Meromorphic functions, Oxford at the Clarendon Press, 1964.
\bibitem{Hay9}
W.K. Hayman, Multivalent functions, 2nd edition, Cambridge Tracts in
Mathematics 110, Cambridge University Press, Cambridge 1994.
\bibitem{HSW}
S. Hellerstein, L.-C. Shen and J. Williamson, 
Real zeros of derivatives of meromorphic functions and solutions of second order differential equations, 
\textit{Trans. Amer. Math. Soc.} 285 (1984), 759-776.
\bibitem{Hil2}
E. Hille, Ordinary differential equations in the complex domain, Wiley,
New York, 1976.
\bibitem{Lai1} I. Laine, Nevanlinna theory and complex differential equations, de Gruyter Studies in Math. 15,  Walter de Gruyter, Berlin/New York 1993.
\bibitem{Lasparse}
J.K. Langley, Bank-Laine functions with sparse zeros, Proc. Amer. Math. Soc. 129 (2001), 1969-1978.
\bibitem{schwarzian} J.K. Langley,  The Schwarzian derivative and the Wiman-Valiron property,  J. Analyse Math.  130 (2016), 71-89.
\bibitem{Lasing2016}
J.K. Langley, Transcendental singularities for a meromorphic function with logarithmic derivative of finite lower order, 
Comput. Methods Funct. Theory 19 (2019), 117-133.
\bibitem{blnewqc}
J.K. Langley, 
Bank-Laine functions, the Liouville transformation and the Eremenko-Lyubich class, to appear, J. Analyse Math.
\bibitem{LV}
O. Lehto and K. Virtanen, Quasiconformal mappings in the plane,
2nd edn., Springer, Berlin, 1973.
\bibitem{Nev2}
R. Nevanlinna, \"Uber Riemannsche Fl\"achen mit endlich vielen Windungspunkten,
Acta Math. 58 (1932), 295-373.
\bibitem{Nev}
R. Nevanlinna, Eindeutige analytische Funktionen,
2. Auflage, Springer, Berlin, 1953.
\bibitem{Ros}
J. Rossi, Second order differential equations with transcendental coefficients,
Proc. Amer. Math. Soc. 97 (1986), 61-66.
\bibitem{Shen} L.C. Shen, Solution to a problem of S. Bank regarding the exponent of convergence of the solutions of a 
differential equation $f''	  + A f= 0$, Kexue Tongbao 30 (1985), 1581-1585.
\bibitem{Shen2}
L.C. Shen, Construction of a differential equation $y'' + Ay = 0$
with solutions having prescribed zeros, Proc. Amer. Math.
Soc. 95 (1985), 544-546.
\bibitem{sixsmithEL} D.J. Sixsmith, Dynamics in the Eremenko-Lyubich class, Conform. Geom. Dyn. 22 (2018), 185-224. 
\bibitem{Steiradial}
N. Steinmetz,
Linear differential equations with exceptional fundamental sets II,
Proc. Amer. Math. Soc. 117 (1993), no. 2, 355--358. 
\bibitem{Tsuji}
M. Tsuji, Potential theory in modern function theory,
Maruzen, Tokyo, 1959.

\end{thebibliography}
}
\end{document}